\documentclass[a4paper,leqno,12pt]{amsproc}
\usepackage{amsfonts}
\usepackage{palatino}
\usepackage{a4wide}
\usepackage{fullpage}
\usepackage{mathpple}

\theoremstyle{plain}
\begingroup

\newtheorem{teo}{Theorem}[section]
\newtheorem{lemma}[teo]{Lemma}
\newtheorem{prop}[teo]{Proposition}
\newtheorem{cor}[teo]{Corollary}
\endgroup

\theoremstyle{definition}
\newtheorem{dfnz}[teo]{Definition}
          
\newtheorem{prob}[teo]{Problem}

\theoremstyle{remark}

\newtheorem{rem}[teo]{Remark}

\newtheorem{ackn}{Acknowledgments\!}

\numberwithin{equation}{section}

\def\div{\operatornamewithlimits{div}\nolimits}
\def\tr{\operatornamewithlimits{tr}\nolimits}
\def\R{{{\mathbb R}}}
\def\SS{{{\mathbb S}}}

\def\ZZ{{{\mathbb Z}}}

\def\RRR{{\mathrm R}}
\def\WWW{{\mathrm W}}
\def\SSS{{\mathrm S}}

\begin{document}

\title[Ricci Solitons -- The Equation Point of View]{Ricci
  Solitons -- The Equation Point of View}

\author[Manolo Eminenti]{Manolo Eminenti}
\address[Manolo Eminenti]{}
\email[M. Eminenti]{manolo@linuz.sns.it}
\author[Gabriele La Nave]{Gabriele La Nave}
\address[Gabriele La Nave]{}
\email[G. La Nave]{lanave@courant.nyu.edu}
\author[Carlo Mantegazza]{Carlo Mantegazza}
\address[Carlo Mantegazza]{Scuola Normale Superiore Pisa, Italy, 56126}
\email[C. Mantegazza]{mantegazza@sns.it}

\date{\today}

\begin{abstract}
We discuss some classification results for Ricci solitons, that
is, self similar solutions of the Ricci Flow.\\
New simpler proofs of some known results will be presented.\\
In detail, we will take the equation point of view, trying to avoid 
the tools provided by considering the dynamic properties of the Ricci
flow.
\end{abstract}

\maketitle

\tableofcontents

\section{Introduction}

\begin{dfnz}
A {\em Ricci Soliton} is a smooth $n$--dimensional complete Riemannian
manifold $(M,g)$ such that there exists a smooth 1--form $\omega$ such that 
\begin{equation}\label{soliton}
2\RRR_{ij}+\nabla_i\omega_j+\nabla_j\omega_i=\frac{2\mu}{n}g_{ij}
\end{equation}
for some constant $\mu\in\R$.

A {\em gradient Ricci Soliton} is a smooth 
$n$--dimensional complete Riemannian
manifold $(M,g)$ such that there exists a smooth function 
$f:M\to\R$, sometimes called {\em potential function}, satisfying
\begin{equation}\label{gsoliton}
\RRR_{ij}+\nabla^2_{ij}f=\frac{\mu}{n}g_{ij}
\end{equation}
for some constant $\mu\in\R$.

Sometimes in literature these manifolds are called {\em quasi--Einstein} manifolds.

A soliton is said to be {\em contracting}, {\em steady} or {\em expanding} if the
constant $\mu\in\R$ is respectively positive, zero or negative.

We say that a soliton is {\em trivial} 
if the form $\omega$ can be chosen to be zero, or the function $f$ to
be constant.\\ 
This is like to say that $(M,g)$ is an Einstein manifold (in
dimension three or if the Weyl tensor is zero, 
it is equivalent to constant curvature).
\end{dfnz}

\begin{rem}
Clearly, a Ricci soliton is a gradient soliton if 
the form $\omega$ above is exact.
\end{rem}

Ricci solitons move under the Ricci flow simply by diffeomorphisms and
homotheties of the initial metric, that is, in other words, they are
stationary points of the Ricci flow in the space of the metrics on $M$
modulo diffeomorphism and scaling.\\
Moreover, their importance is due to the fact that they arise as
blow--up limits of the Ricci flow when singularities develop.\\
In this note we try to analyze these manifolds
forgetting the properties of the Ricci flow and the related 
{\em dynamical} techniques, just looking at the defining elliptic
equations~\eqref{soliton} and~\eqref{gsoliton}.

\smallskip

We suggest to the interested reader 
the preprint of Derdzinski~\cite{derdz2} and the survey of
Cao~\cite{cao3} as very comprehensive reviews of the present
literature, open problems and recent developments.

\smallskip

During the publication of this paper, there appear 
several preprints extending the results (or presenting different proofs)
to the complete noncompact case~\cite{caowang,fmanzan,na1,nw1,nw2,pw}.

\begin{ackn} We wish to thank Xiaodong Cao for pointing us at some 
inaccuracies in earlier versions of the paper.
\end{ackn}

\section{General Computations}

We recall some well known facts from Riemannian geometry:
\begin{itemize}
\item The Schur's Lemma, when $n>2$
$$
2{\mathrm {div}}{\mathrm {Ric}}=d\RRR\,.
$$
\item The fact that $\RRR_{ij}=g_{ij}\RRR/2$ when $n=2$.
\item The formula for the interchange of covariant derivatives of a
  form,
$$
\nabla^2_{ij}\omega_k-\nabla^2_{ji}\omega_k= \RRR_{ijks}\omega^s\,.
$$
\item The decomposition of the Riemann tensor, 
$$
\RRR_{ijkl}=-\frac{\RRR}{(n-1)(n-2)}(g_{ik}g_{jl}-g_{il}g_{jk}) +
\frac{1}{n-2}(\RRR_{ik}g_{jl}-\RRR_{il}g_{jk}
+\RRR_{jl}g_{ik}-\RRR_{jk}g_{il})+\WWW_{ijkl}\,.
$$
\item The fact that the Weyl tensor $\WWW$ is zero when $n\leq3$.
\end{itemize}

Now we work out some consequences of equation~\eqref{gsoliton}.

\begin{prop} Let $(M,g)$ be a Ricci gradient soliton, then the following formulas hold,
\begin{equation}\label{equ1}
\RRR+\Delta f=\mu
\end{equation}
\begin{equation}\label{equ2}
\nabla_i \RRR=2\RRR_{ij}\nabla^jf
\end{equation}
\begin{equation}\label{commut1}
\nabla_j \RRR_{ik}-\nabla_i\RRR_{jk}=\RRR_{ijks}\nabla^sf\,.
\end{equation}
\begin{equation}\label{equ8}
\RRR+\vert\nabla f\vert^2-\frac{2\mu}{n}f=\,\text{constant}
\end{equation}
\begin{equation}\label{equ4} 
\Delta\RRR=\langle\nabla\RRR\,\vert\,\nabla
f\rangle+\frac{2\mu}{n}\RRR-2\vert{\mathrm {Ric}}\vert^2
\end{equation}
\begin{equation}\label{equ3}
\Delta\RRR_{ij}=\langle\nabla\RRR_{ij}\,\vert\,\nabla
f\rangle+\frac{2\mu}{n}\RRR_{ij}-2\RRR_{ikjs}\RRR^{ks}
\end{equation}
\begin{align}
\Delta\RRR_{ik}
=&\,\langle\nabla\RRR_{ik}\,\vert\,\nabla
f\rangle+\frac{2\mu}{n}\RRR_{ik}-\WWW_{ijkl}R^{jl}\label{equ6}\\
&\,+\frac{2}{(n-1)(n-2)}
\bigl(\RRR^2 g_{ik}-n\RRR\RRR_{ik}
+2(n-1)\SSS_{ik}-(n-1)\SSS g_{ik}\bigr)\,.\nonumber
\end{align}
where $\SSS_{ik}=\RRR_{ij}g^{jl}\RRR_{lk}$ and
$\SSS=\tr\SSS=\vert{\mathrm {Ric}}\vert^2$.
\end{prop}
\begin{proof}
Equation~\eqref{equ1}: we simply contract equation~\eqref{gsoliton}.

\smallskip

\noindent Equation~\eqref{equ2}: we take the divergence of the Ricci tensor, by using
equation~\eqref{gsoliton},
\begin{align*}
\div{\mathrm{Ric}}_i=&\,g^{jk}\nabla_k\RRR_{ij}\\
=&\,-g^{jk}\nabla_k\nabla_i\nabla_jf\\
=&\,-g^{jk}\nabla_i\nabla_k\nabla_jf-g^{jk}\RRR_{kijs}\nabla^sf\\
=&\,-\nabla_i\Delta f-\RRR_{is}\nabla^sf\,,
\end{align*}
where we used the formula for the interchange of covariant
derivatives.
Using then equation~\eqref{equ1} and Schur's Lemma we get
\begin{equation*}
\frac{1}{2}\nabla_i\RRR=-\nabla_i(\mu-\RRR) -\RRR_{is}\nabla^sf=
\nabla_i \RRR-\RRR_{is}\nabla^sf\,,
\end{equation*}
hence, equation~\eqref{equ2} follows.

\smallskip

\noindent Equation~\eqref{commut1}: it follows by a computation analogous
to the previous one.

\smallskip

\noindent Relation~\eqref{equ8}: it follows by simply differentiating
the quantity on the left side and using equations~\eqref{equ2}
and~\eqref{gsoliton}.

\smallskip

\noindent Equation~\eqref{equ4}: once obtained equation~\eqref{equ3},
it follows by contracting it with the metric $g$.

\smallskip

\noindent Equation~\eqref{equ3}: we compute the Laplacian of the Ricci
tensor by means of equation~\eqref{commut1} and the second Bianchi
identity,
\begin{align*}
\Delta\RRR_{ik}=\nabla^j\nabla_j\RRR_{ik}
=&\,\nabla^j\nabla_i\RRR_{jk}+\nabla^j\RRR_{ijks}\nabla^sf+\RRR_{ijks}\nabla^{j}\nabla^sf\\
=&\,\nabla_i\nabla^j\RRR_{jk}+\RRR^j_{\phantom{j}ijs}\RRR^s_k+\RRR^j_{\phantom{j}iks}\RRR^s_j\\
&\,+\nabla_k\RRR^j_{\phantom{j}sij}\nabla^sf-\nabla_s\RRR^j_{\phantom{j}kij}\nabla^sf
+\RRR_{ijks}\nabla^{j}\nabla^sf\\
=&\,\nabla_i\nabla^j\RRR_{jk}+\RRR_{is}\RRR^s_k+\RRR^j_{\phantom{j}iks}\RRR^s_j\\
&\,-\nabla_k\RRR_{si}\nabla^sf+\nabla_s\RRR_{ki}\nabla^sf
+\RRR_{ijks}\nabla^{j}\nabla^{s}f\\
=&\,\frac{1}{2}\nabla_i\nabla_k\RRR+\RRR_{is}\RRR^s_k+\RRR^j_{\phantom{j}iks}\RRR^s_j\\
&\,-\nabla_k\RRR_{si}\nabla^sf+\langle\nabla\RRR_{ik}\,\vert\,\nabla f\rangle
-\RRR_{ijks}\RRR^{js}+\frac{\mu}{n}\RRR_{ik}\,,\\
=&\,\langle\nabla\RRR_{ik}\,\vert\,\nabla
f\rangle+\frac{\mu}{n}\RRR_{ik}-2\RRR_{ijks}\RRR^{js}
+\RRR_{is}\RRR^s_k\\
&\,+\frac{1}{2}\nabla_k\nabla_i\RRR-\nabla_k\RRR_{is}\nabla^sf\,,\\
\end{align*}
where we used also Schur's Lemma, substituted $\nabla^j\nabla^sf$ by means of
equation~\eqref{gsoliton} and rearranged some terms in the last line.\\
Differentiating relation~\eqref{equ2} we obtain
$$
0=\frac{1}{2}\nabla_k\left(\nabla_i\RRR-2\RRR_{is}\nabla^sf\right)
$$
hence,
$$
\frac{1}{2}\nabla_k\nabla_i\RRR-\nabla_k\RRR_{is}\nabla^sf=\RRR_{is}\nabla^s\nabla_kf\,.
$$
Then, we conclude
\begin{align*}
\Delta\RRR_{ik}
=&\,\langle\nabla\RRR_{ik}\,\vert\,\nabla
f\rangle+\frac{\mu}{n}\RRR_{ik}-2\RRR_{ijks}\RRR^{js}
+\RRR_{is}\RRR^s_k+\RRR_{is}\nabla^s\nabla_kf\\
=&\,\langle\nabla\RRR_{ik}\,\vert\,\nabla
f\rangle+\frac{\mu}{n}\RRR_{ik}-2\RRR_{ijks}\RRR^{js}
+\RRR_{is}\RRR^s_k+\frac{\mu}{n}\RRR_{is}-\RRR_{is}\RRR^s_k\\
=&\,\langle\nabla\RRR_{ik}\,\vert\,\nabla
f\rangle+\frac{2\mu}{n}\RRR_{ik}-2\RRR_{ijks}\RRR^{js}\,.
\end{align*}

\smallskip

\noindent Equation~\eqref{equ6}: by using the decomposition of the
Riemann tensor we can go on with the computation in
formula~\eqref{equ3},
\begin{align*}
\Delta\RRR_{ik}
=&\,\langle\nabla\RRR_{ik}\,\vert\,\nabla
f\rangle+\frac{2\mu}{n}\RRR_{ik}
-2\WWW_{ijkl}\RRR^{jl}\\
&\,+\frac{2\RRR}{(n-1)(n-2)}(g_{ik}g_{jl}-g_{il}g_{jk})\RRR^{jl}\\ 
&\,-\frac{2}{n-2}(\RRR_{ik}g_{jl}-\RRR_{il}g_{jk}+\RRR_{jl}g_{ik}
-\RRR_{jk}g_{il})\RRR^{jl}\\
=&\,\langle\nabla\RRR_{ik}\,\vert\,\nabla
f\rangle+\frac{2\mu}{n}\RRR_{ik}
-2\WWW_{ijkl}\RRR^{jl}\\
&\,+\frac{2\RRR}{(n-1)(n-2)}(\RRR g_{ik}-\RRR_{ik})
-\frac{2}{n-2}(\RRR\RRR_{ik}-2\SSS_{ik}+\SSS g_{ik})\\
=&\,\langle\nabla\RRR_{ik}\,\vert\,\nabla
f\rangle+\frac{2\mu}{n}\RRR_{ik}
-2\WWW_{ijkl}\RRR^{jl}\\
&\,+\frac{2}{(n-1)(n-2)}
\bigl(\RRR^2 g_{ik}-n\RRR\RRR_{ik}
+2(n-1)\SSS_{ik}-(n-1)\SSS g_{ik}\bigr)\,.
\end{align*}
\end{proof}

\subsection{The Case $n=2$}

We can use the complete description of the curvature tensor via the
scalar curvature $\RRR$, that is $\RRR_{ij}=\RRR g_{ij}/2$.\\

\begin{prop} If $n=2$ there hold,
\begin{align*}
\nabla^2_{ij} f=&\,\frac{\mu-\RRR}{2}g_{ij}\\
\nabla\RRR=&\,\RRR\nabla f\\
\Delta\RRR=&\,\langle\nabla\RRR\,\vert\,\nabla f\rangle+\mu\RRR-\RRR^2
\end{align*}
\end{prop}

\subsection{The Case $n=3$}

Using equation~\eqref{equ6}, as the Weyl tensor $\WWW$ is 
identically zero for every 3--manifold, we get
\begin{prop} If $n=3$ there hold,
\begin{equation*}
\Delta\RRR_{ik}
= \langle\nabla\RRR_{ik}\,\vert\,\nabla
f\rangle+\frac{2\mu}{3}\RRR_{ik}+\RRR^2 g_{ik}-3\RRR\RRR_{ik}
+4\SSS_{ik}-2\SSS g_{ik}
\end{equation*}
with $\SSS_{ik}=\RRR_{ij}g^{jl}\RRR_{lk}$ and
$\SSS=\vert{\mathrm {Ric}}\vert^2$.
\end{prop}
This proposition clearly generalizes to dimension $n>3$
  when $\WWW=0$,
\begin{align}
\Delta\RRR_{ik}
=&\,\langle\nabla\RRR_{ik}\,\vert\,\nabla
f\rangle+\frac{2\mu}{n}\RRR_{ik}\label{equ5}\\
&\,+\frac{2}{(n-1)(n-2)}
\bigl(\RRR^2 g_{ik}-n\RRR\RRR_{ik}
+2(n-1)\SSS_{ik}-(n-1)\SSS g_{ik}\bigr)\,.\nonumber
\end{align}

\section{Compact Ricci Solitons}

By means of Perelman work~\cite{perel1} and previous others, see
Hamilton~\cite{hamilton5} (dimension two) and Ivey~\cite{ivey1}
(dimension 3), we have the following fact.

\begin{teo}[Perelman]\label{p4} 
Every compact Ricci soliton is a gradient Ricci soliton.
\end{teo}
\begin{proof}
Let $(M,g)$ be a Ricci Soliton, with a potential form $\omega$. 
We start with the following computation for a generic smooth function $f:M\to\R$,
\begin{align*}
g^{kj}\nabla_k&\,[2(\RRR_{ij}+\nabla^2_{ij}f-\mu g_{ij}/n)e^{-f}]\\
=&\,(\nabla_i\RRR+2\Delta\nabla_if)e^{-f}
-2[(\RRR_{ij}+\nabla^2_{ij}f-\mu g_{ij}/n)g^{jk}\nabla_k f]e^{-f}\\
=&\,(\nabla_i\RRR+2\nabla_i\Delta f +2\RRR_{is}\nabla^s f)e^{-f}
-2[(\RRR_{ij}+\nabla^2_{ij}f-\mu g_{ij}/n)g^{jk}\nabla_k f]e^{-f}\\
=&\,(\nabla_i\RRR+2\nabla_i\Delta f-2g^{jk}\nabla^2_{ij}f\nabla_k f
+2\mu/n\nabla_i f) e^{-f}\\
=&\,\nabla_i(\RRR+2\Delta f-\vert\nabla f\vert^2+2\mu f/n) e^{-f}\,.
\end{align*}
Hence, supposing to find a smooth function $f:M\to\R$ such that
\begin{equation}\label{equ100}
\RRR+2\Delta f-\vert\nabla f\vert^2+2\mu f/n
\end{equation}
is constant, we have
$$
\div[(\mathrm{Ric}+\nabla^2 f-\mu g/n)e^{-f}]=0\,.
$$
Then, as $\nabla_l\omega_k+\nabla_k\omega_l=-2\RRR_{lk}+2\mu g_{lk}/n$,
\begin{align*}
\div[(\nabla_k&\,
f-\omega_k)g^{kj}(\RRR_{ij}+\nabla^2_{ij}f-\mu g_{ij}/n)e^{-f}]\\
=&\,(\nabla^2_{lk}f-\nabla_l\omega_k)g^{kj}g^{li}
(\RRR_{ij}+\nabla^2_{ij}f-\mu g_{ij}/n)e^{-f}\\
=&\,(2\nabla^2_{lk}f-\nabla_l\omega_k-\nabla_k\omega_l)g^{kj}g^{li}
(\RRR_{ij}+\nabla^2_{ij}f-\mu g_{ij}/n)e^{-f}/2\\
=&\,(2\nabla^2_{lk}f+2\RRR_{lk}-2\mu g_{lk}/n)g^{kj}g^{li}
(\RRR_{ij}+\nabla^2_{ij}f-\mu g_{ij}/n)e^{-f}/2\\
=&\,\vert\RRR_{ij}+\nabla^2_{ij}f-\mu g_{ij}/n\vert^2 e^{-f}\,,
\end{align*}
where, passing from the second to the third line, we substituted
$\nabla_l\omega_k$ with $(\nabla_l\omega_k+\nabla_k\omega_l)/2$ since
the skew--symmetric component of $\nabla\omega$ vanishes once we
contract it with the symmetric 2--form $\RRR_{ij}+\nabla^2_{ij}f-\mu
g_{ij}/n$.\\
Hence, we conclude that
$$
0\leq Q=\vert\mathrm{Ric}+\nabla^2 f-\mu g/n\vert^2
e^{-f}=\div\mathrm{T}
$$
for some 1--form $\mathrm{T}$.\\
Integrating $Q$ on $M$, we immediately get that $Q=0$, since it is 
nonnegative.\\
This clearly implies that $(M,g)$ is a gradient Ricci soliton with a
potential $f$.

The existence of a function $f$ such that relation~\eqref{equ100}
holds, can be proven by constrained minimization of Perelman's
${\mathcal{W}}$ functional (defined in~\cite{perel1}). A logarithmic
Sobolev estimate is needed, see Appendix~\ref{appA}.
\end{proof}

\begin{prob} Is it possible to prove Theorem~\ref{p4}
  without using arguments related to Ricci flow, that is, 
showing directly that the form  $\omega$ in equation~\eqref{soliton}
is exact?
\end{prob}

\begin{rem} In the noncompact case, there exists {\em non} gradient
 Ricci solitons, see Baird and Danielo~\cite{bairdan}, Lott~\cite{lott}.
\end{rem}

We can then concentrate ourselves on compact {\em gradient} Ricci
solitons.\\
The key tool will be the maximum principle for elliptic equations.

We start from equation~\eqref{equ4},
\begin{equation*}
\Delta\RRR=\langle\nabla\RRR\,\vert\,\nabla
f\rangle+\frac{2\mu}{n}\RRR-2\vert{\mathrm {Ric}}\vert^2\,
\end{equation*}
noticing that 
$\mu=\frac{1}{{\mathrm {Vol}}(M)}\int_M\RRR\geq\RRR_{\min}$, with
equality if and only if $\RRR$ is constant.\\
We have,
\begin{align*}
\Delta\RRR=&\,\langle\nabla\RRR\,\vert\,\nabla
f\rangle+\frac{2\mu}{n}\RRR
-2\vert{{\overset{\circ}{\mathrm {Ric}}}}\vert^2-2\RRR^2/n\\
\leq&\,\langle\nabla\RRR\,\vert\,\nabla
f\rangle+\frac{2\RRR}{n}(\mu-\RRR)
\end{align*}
where we denoted with ${\overset{\circ}{\mathrm {Ric}}}$ the {\em tracefree} part of the Ricci tensor.\\
When $\RRR$ gets its minimum,
\begin{equation*}
\Delta\RRR_{\min}\leq\frac{2\RRR_{\min}}{n}(\mu-\RRR_{\min})\,.
\end{equation*}
This relation, by the strong maximum principle, implies that if $\RRR$
is nonconstant, then it must be positive everywhere, hence also $\mu$ is
positive.

\begin{prop} Every steady or expanding compact Ricci soliton has 
the scalar curvature $\RRR$ constant (and equal to the constant $\mu$).
\end{prop}
 
Coming back to equation~\eqref{equ1} this fact forces 
$\Delta f=0$, hence, since we are on a compact manifold, $f$ is
constant and the soliton is trivial.

\begin{cor}\label{p1} 
Every steady or expanding compact Ricci soliton
  is trivial.
\end{cor}

Now we deal with contracting compact gradient Ricci solitons.\\

The two--dimensional case is special, we saw that $\RRR>0$, hence
topologically we are dealing with $\SS^2$ or its $\ZZ_2$--quotient
${\mathbb{RP}^2}$.\\
The relevant equation is
$$
\nabla^2_{ij} f=\frac{\mu-\RRR}{2}g_{ij}\,,
$$
indeed, differently by the 3--dimensional case, the II Bianchi
identity (hence the Schur's Lemma) is a void condition, making this
case more difficult.

The following result was first proved by Hamilton~\cite{hamilton5} with
an argument using the uniformization theorem which can be strongly 
simplified by means of Kazdan--Warner identity (which relies on
uniformization), see~\cite[pag.~131]{chknopf} and~\cite{choyau}. 
Recently, Chen, Lu and Tian found a simple proof independent by
uniformization of surfaces~\cite{chenlutian}.

\begin{prop}\label{p3} 
Every contracting, compact, two--dimensional 
Ricci soliton is ${\mathbb S}^2$ or ${\mathbb{RP}^2}$ 
with the standard metric.
\end{prop}

We assume now to be in dimension three or in higher dimension with a
zero Weyl tensor. As we said, the scalar curvature 
$\RRR$ must be positive everywhere, then by means of
equation~\eqref{equ5} we have,
\begin{align*}
\Delta\left(\frac{\RRR_{ik}}{\RRR}\right)
=&\,\frac{\Delta\RRR_{ik}}{\RRR}-\frac{\RRR_{ik}\Delta\RRR}{\RRR^2}
+2\frac{\vert\nabla\RRR\vert^2\RRR_{ik}}{\RRR^3}
-2\frac{\langle\nabla\RRR_{ik}\,\vert\,\nabla\RRR\rangle}{\RRR^2}\,,
\end{align*}
substituting the equations for $\Delta\RRR_{ik}$ and $\Delta\RRR$ we get
\begin{align*}
\Delta\left(\frac{\RRR_{ik}}{\RRR}\right)
=&\,\frac{\Delta\RRR_{ik}}{\RRR}-\frac{\RRR_{ik}\Delta\RRR}{\RRR^2}
+2\frac{\vert\nabla\RRR\vert^2\RRR_{ik}}{\RRR^3}
-2\frac{\langle\nabla\RRR_{ik}\,\vert\,\nabla\RRR\rangle}{\RRR^2}\\
=&\,\frac{\langle\nabla\RRR_{ik}\,\vert\,\nabla
f\rangle}{\RRR}+\frac{2\mu}{n}\frac{\RRR_{ik}}{\RRR}\\
&\,+\frac{2\bigl(\RRR^2 g_{ik}-n\RRR\RRR_{ik}
+2(n-1)\SSS_{ik}-(n-1)\SSS g_{ik}\bigr)}{(n-1)(n-2)\RRR}\\
&\,-\frac{\RRR_{ik}\Delta\RRR}{\RRR^2}
+2\frac{\vert\nabla\RRR\vert^2\RRR_{ik}}{\RRR^3}
-2\frac{\langle\nabla\RRR_{ik}\,\vert\,\nabla\RRR\rangle}{\RRR^2}\\
=&\,\Bigl\langle\frac{\nabla{\RRR_{ik}}}{\RRR}
\,\Bigr\vert\,\nabla f -2\frac{\nabla\RRR}{\RRR}\Bigr\rangle
+\frac{2\mu}{n}\frac{\RRR_{ik}}{\RRR}\\
&\,+\frac{2\bigl(\RRR^2 g_{ik}-n\RRR\RRR_{ik}
+2(n-1)\SSS_{ik}-(n-1)\SSS g_{ik}\bigr)}{(n-1)(n-2)\RRR}\\
&\,-\frac{\RRR_{ik}\langle\nabla\RRR\,\vert\,\nabla
f\rangle+\frac{2\mu}{n}\RRR\RRR_{ik}-2\vert{\mathrm
  {Ric}}\vert^2\RRR_{ik}}{\RRR^2}
+2\frac{\vert\nabla\RRR\vert^2\RRR_{ik}}{\RRR^3}\\
=&\,\Bigl\langle\nabla\Bigl(\frac{\RRR_{ik}}{\RRR}\Bigr)
\,\Bigr\vert\,\nabla f -2\frac{\nabla\RRR}{\RRR}\Bigr\rangle\\
&\,+\frac{2}{(n-1)(n-2)}\left(\RRR g_{ik}-n\RRR_{ik}
+\frac{2(n-1)\SSS_{ik}}{\RRR}-\frac{(n-1)\SSS g_{ik}}{\RRR}\right)
+\frac{2\SSS\RRR_{ik}}{\RRR^2}\\
=&\,\Bigl\langle\nabla\Bigl(\frac{\RRR_{ik}}{\RRR}\Bigr)
\,\Bigr\vert\,\nabla f -2\frac{\nabla\RRR}{\RRR}\Bigr\rangle\\
&\,+\frac{2\left(\RRR^3g_{ik}
-n\RRR^2\RRR_{ik}
+2(n-1)\RRR\SSS_{ik}
-(n-1)\SSS\RRR g_{ik}+(n-1)(n-2)\SSS\RRR_{ik}\right)}{(n-1)(n-2)\RRR^2}\,.
\end{align*}
Let now $\lambda_{\min}:M\to\R$ to be 
the minimal eigenvalue of the Ricci tensor. If $p\in M$ is 
the point where $\lambda_{\min}/\RRR$ gets its minimum with
eigenvector $v_p$, we consider a
local unit smooth tangent vector field $w=w^i$ such that $w(p)=v_p$,
$\nabla w^i(p)=\Delta w^i(p)=0$. Then the smooth function
$\RRR_{ij}w^iw^j/\RRR$ has a local minimum at $p$, 
$(\RRR_{ij}w^iw^j/\RRR)(p)=\lambda_{\min}(p)/\RRR(p)$,
$\nabla(\RRR_{ij}w^iw^j/\RRR)(p)=0$ and 
$\Delta(\RRR_{ij}w^iw^j/\RRR)(p)\geq 0$.\\
By the assumptions on the derivatives of $w$ at $p$ we have$\nabla(\RRR_{ij}/\RRR)(p)v_p^iv_p^j=0$ and 
$\Delta(\RRR_{ij}/\RRR)(p)v_p^iv_p^j\geq 0$, hence, using the previous
equation,
$$
0\leq \Delta(\RRR_{ij}/\RRR)(p)v_p^iv_p^j=\frac{2\left(\RRR^3-n\lambda_{\min}\RRR^2+2(n-1)\lambda_{\min}^2\RRR-(n-1)\SSS\RRR
+(n-1)(n-2)\lambda_{\min}\SSS\right)}{\RRR^2}\,,
$$
where the right hand side is evaluated at $p\in M$.\\
This implies,
\begin{equation}\label{equ9}
0\leq 
\RRR^3-n\lambda_{\min}\RRR^2+2(n-1)\lambda_{\min}^2\RRR-(n-1)\SSS\RRR
+(n-1)(n-2)\lambda_{\min}\SSS\,.
\end{equation}
We work on this term, setting $\widetilde{\RRR}$ and 
$\widetilde{\SSS}$ to be respectively the sum and the sum of the 
squares of the eigenvalues of the Ricci tensor, but $\lambda_{\min}$.
\begin{align*}
0\leq&\, \RRR^3-n\lambda_{\min}\RRR^2+2(n-1)\lambda_{\min}^2\RRR-(n-1)\SSS\RRR
+(n-1)(n-2)\lambda_{\min}\SSS\\
=&\,(\lambda_{\min}+\widetilde{\RRR})^3
-n\lambda_{\min}(\lambda_{\min}+\widetilde{\RRR})^2
+2(n-1)\lambda_{\min}^2(\lambda_{\min}+\widetilde{\RRR})\\
&\,-(n-1)(\lambda_{\min}^2+\widetilde{\SSS})(\lambda_{\min}+\widetilde{\RRR})
+(n-1)(n-2)\lambda_{\min}(\lambda_{\min}^2+\widetilde{\SSS})\\
=&\,\lambda_{\min}^3\bigl(1-n+2(n-1)-(n-1)+(n-1)(n-2)\bigr)\\
&\,+\lambda_{\min}^2\bigl(3\widetilde{\RRR}-2n\widetilde{\RRR}
+2(n-1)\widetilde{\RRR}-(n-1)\widetilde{\RRR}\bigr)\\
&\,+\lambda_{\min}\bigl(3\widetilde{\RRR}^2-n\widetilde{\RRR}^2
-(n-1)\widetilde{\SSS}+(n-1)(n-2)\widetilde{\SSS}\bigr)\\
&\,+(\widetilde{\RRR}^3-(n-1)\widetilde{\RRR}\widetilde{\SSS}\bigr)\\
=&\,(n-1)(n-2)\lambda_{\min}^3-(n-2)\lambda_{\min}^2\widetilde{\RRR}\\
&\,+(n-3)\lambda_{\min}\bigl((n-1)\widetilde{\SSS}-\widetilde{\RRR}^2\bigr)
-\widetilde{\RRR}\bigl((n-1)\widetilde{\SSS}-\widetilde{\RRR}^2\bigr)\\
=&\,(n-2)\lambda_{\min}^2\bigl(\lambda_{\min}(n-1)-\widetilde{\RRR}\bigr)
+\bigl((n-3)\lambda_{\min}-\widetilde{\RRR}\bigr)
\bigl((n-1)\widetilde{\SSS}-\widetilde{\RRR}^2\bigr)\,.
\end{align*}
Now, as $\RRR$ is positive, both terms 
$\bigl(\lambda_{\min}(n-1)-\widetilde{\RRR}\bigr)$ and 
$\bigl((n-3)\lambda_{\min}-\widetilde{\RRR}\bigr)$ are nonpositive. 
Using the {\em Arithmetic--Quadratic} mean inequality we see that 
the term $\bigl((n-1)\widetilde{\SSS}-\widetilde{\RRR}^2\bigr)$ must
be nonnegative, so we conclude that all this expression is
nonpositive. Hence, it must be zero at $p$, the point where
$\lambda_{\min}/\RRR$ gets its minimum.\\
There are only two possibilities this can happen: 
either $\lambda_{\min}(p)=0$ and all the other $n-1$ eigenvalues of the 
Ricci tensor are equal, or all the eigenvalues are equal (and positive as $R>0$).\\ 
In this latter case, $\lambda_{\min}(p)=\RRR/n$ and since we are in
the point of minimum, $\RRR_{ij}/\RRR\geq g_{ij}/n$ or $\RRR_{ij}\geq
\RRR g_{ij}/n$ on the whole manifold. But this inequality easily implies
that $(M,g)$ is an 
Einstein manifold (the soliton is trivial).\\
In the other case, as $\lambda_{\min}(p)=0$ and all the other $n-1$ eigenvalues of the 
Ricci tensor are equal to $\RRR/(n-1)$. It can be shown that locally around $p$ 
the eigenvector $v(q)$ realizing the minimal eigenvalue
$\lambda_{\min}(q)$ can be chosen smoothly depending on the point $q$
(locally there are no ``bifurcations'' of the minimal eigenvalue of
$\RRR_{ij}$). Then, as 
$\RRR_{ij}v^i=\lambda_{\min}v_j$, differentiating this relation, we get
$$
\nabla^k\RRR_{ij}v^i
=\nabla^k\lambda_{\min}v_j+\lambda_{\min}\nabla^kv_j-\RRR_{ij}\nabla^kv^i\,,
$$
then we compute for $\lambda_{\min}=\RRR_{ij}v^iv^j$,
\begin{align*}
\Delta\lambda_{\min}=&\,\Delta(\RRR_{ij}v^iv^j)\\
=&\,\Delta\RRR_{ij}v^iv^j+4\nabla^k\RRR_{ij}v^i\nabla_kv^j
+2\RRR_{ij}\nabla^kv^i\nabla_kv^j+2\RRR_{ij}v^i\Delta v^j\\
=&\,\Delta\RRR_{ij}v^iv^j
+2\RRR_{ij}\nabla^kv^i\nabla_kv^j+2\RRR_{ij}v^i\Delta v^j\\
&\,+4\nabla^k\lambda_{\min}v_j\nabla_kv^j
+4\lambda_{\min}\nabla^kv_j\nabla_kv^j-4\RRR_{ij}\nabla^kv^i\nabla_kv^j\\
=&\,\Delta\RRR_{ij}v^iv^j
-2\RRR_{ij}\nabla^kv^i\nabla_kv^j+2\RRR_{ij}v^i\Delta v^j
+2\nabla^k\lambda_{\min}\nabla_k\vert v\vert^2
+4\lambda_{\min}\vert\nabla v\vert^2\\
\leq&\,\Delta\RRR_{ij}v^iv^j
-2\lambda_{\min}\nabla^kv_j\nabla^kv^j
+2\RRR_{ij}v^i\Delta v^j
+4\lambda_{\min}\vert\nabla v\vert^2\\
=&\,\Delta\RRR_{ij}v^iv^j
+2\lambda_{\min}\vert\nabla v\vert^2
+2\lambda_{\min}v_j\Delta v^j\\
=&\,\Delta\RRR_{ij}v^iv^j
+\lambda_{\min}\Delta\vert v\vert^2\\
=&\,\Delta\RRR_{ij}v^iv^j\,.
\end{align*}
Working as before, we obtain the following elliptic inequality for the minimal
eigenvalue, locally around $p$,
\begin{align*}
\Delta\lambda_{\min}\leq &\,\langle\nabla\lambda_{\min}\,\vert\,\nabla f\rangle
+\frac{2\mu}{n}\lambda_{\min}\\
&\,+\frac{2}{(n-1)(n-2)}
\bigl(\RRR^2-n\RRR\lambda_{\min}
+2(n-1)\lambda_{\min}^2-(n-1)\SSS\bigr)\,.
\end{align*}
This inequality implies
\begin{align*}
\Delta\Bigl(\frac{\lambda_{\min}}{\RRR}\Bigr)
\leq&\,\Bigl\langle\nabla\Bigl(\frac{\lambda_{\min}}{\RRR}\Bigr)
\,\Bigr\vert\,\nabla f -2\frac{\nabla\RRR}{\RRR}\Bigr\rangle\\
&\,+\frac{2}{(n-1)(n-2)}
\bigl(\RRR-n\lambda_{\min}
+2(n-1)\lambda_{\min}^2/\RRR-(n-1)\SSS/\RRR\bigr)+2\lambda_{\min}\SSS/\RRR^2\,,
\end{align*}
holding locally around $p$ where $\lambda_{\min}/\RRR$ get the local
minimum zero.\\
As at this local minimum
$\Delta\Bigl(\frac{\lambda_{\min}}{\RRR}\Bigr)=0$, by 
the strong maximum principle, it follows that $\lambda_{\min}/\RRR$ is
locally constant around $p$. Then it is an easy consequence that this
must hold on the whole connected $M$ and $\lambda_{\min}$ is
identically zero.\\
Getting back to the initial equation~\eqref{gsoliton}, we put 
ourselves at the point $p\in M$ where the function 
$f$ gets its maximum. If $v\in T_pM$ is the unit zero eigenvector of
the Ricci tensor, we take normal coordinates at $p$ such that
$v=\partial_{x_1}$, hence
\begin{equation*}
\RRR_{ij}(p)v^iv^j+\nabla^2_{ij}f(p)v^iv^i=\frac{\mu}{n}\,,
\end{equation*}
that is,
\begin{equation*}
\nabla^2_{11}f(p)=\mu/n>0\,,
\end{equation*}
which is impossible as $p$ is a maximum point for $f$.

\begin{prop}\label{p2} 
Every contracting, compact, three--dimensional 
Ricci soliton is a quotient of ${\mathbb S}^3$ with the standard
metric.
\end{prop}

\begin{prop}\label{p5}
Every contracting, compact, Ricci soliton when $n>3$
and the Weyl tensor is zero, is trivial. Then, it is a quotient of
${\mathbb S}^n$ with the standard metric.
\end{prop}
\begin{rem} In the recent preprint~\cite{caowang} Cao and Wang also prove this result by means of a completely different method.
\end{rem}

When $n>3$, we have counterexamples to the triviality of compact Ricci
solitons due to Koiso~\cite{koiso1}, Cao~\cite{cao1}, Feldman, Ilmanen
and Ni~\cite{ilman6}. Moreover, some of these examples have ${\mathrm
  {Ric}}>0$. See also Bryant~\cite{bry1}.\\
In general, we only know that it must be $\RRR>0$ and nonconstant.

\begin{prob}
Are there special conditions in dimension $n=4$ (on the Weyl tensor?)
assuring that a contracting, compact, Ricci soliton is trivial?
\end{prob}

\begin{prob}
Are there counterexamples in dimension $n=5$?
\end{prob}

Recently B\"ohm and Wilking~\cite{bohmwilk} showed the following
Hamilton's conjecture.
\begin{teo} If the Riemann {\em operator} is definite positive,
every contracting, compact, Ricci soliton is trivial.
\end{teo}
Hence, it is a quotient of ${\mathbb S}^n$, by a theorem of
Tachibana~\cite{tachib1}.

Previously, by Hamilton work~\cite{hamilton2} this result was known for
$n\leq 4$ and there was a partial result by Cao~\cite{caox1} in any
dimension.

\begin{prob}
The sectional curvatures of a compact, NONtrivial, contracting Ricci soliton 
can be all positive (nonnegative)?
\end{prob}

\begin{prob}
What are in general the properties of compact, NONtrivial, 
contracting Ricci solitons?
\end{prob}
See this paper by Derdzinski~\cite{derdz1}.\\
We are also aware of a preprint~\cite{fergarciario} of
Fern\'andez--L\'opez and Garc\'{\i}a--R\'{\i}o where they show 
that the first fundamental group has to be finite.\\
We give here a short
proof of this fact.

\begin{prop} The first fundamental group of a compact shrinking Ricci
  solitons is finite.
\end{prop}
\begin{proof}
Denoting with $\pi:\widetilde{M}\to M$ the Riemannian universal
covering of $M$, it is well known that the fundamental group is in
one--to--one correspondence with the discrete counterimage of a
basepoint $p\in M$. Clearly $\widetilde{M}$ is again a shrinking gradient Ricci
soliton (possibly noncompact) with a potential function
$\widetilde{f}=f\circ\pi$, since $\pi$ is a local isometry.

Let $a$ and $b$ be a pair of points with $\pi(a)=\pi(b)=p$ and
$\gamma:[0,L]\to\widetilde{M}$ the minimal geodesic between them, 
parametrized by arclength.\\ 
Let $E_i$ be an orthonormal basis of $T_a\widetilde{M}$ such that
$E_1=\gamma^\prime(0)$ extended by parallel transport along
$\gamma$. If $h(t)=\sin\left({\frac{\pi t}{L}}\right)$, we define the fields
$Y_i(t)=h(t)E_i$, zero at $t=0,L$.\\
As $\gamma$ is minimal, the index form is nonnegative definite,
\begin{equation*}
0\leq I(Y_i,Y_i)=\int_0^L \vert
Y_i^\prime\vert^2-\RRR(Y_i,\gamma',Y_i,\gamma')\,dt=
\int_0^L \vert
h'(t)\vert^2-h^2(t)\RRR(E_i,\gamma',E_i,\gamma')\,dt\,,
\end{equation*}
and after summing on $i\in{1,\dots,n}$, we get
\begin{equation*}
0\leq\frac{\pi^2(n-1)}{L^2}\int_0^L \cos^2\left({\frac{\pi t}{L}}\right)\,dt 
-\int_0^L \sin^2\left({\frac{\pi t}{L}}\right)\mathrm{Ric}(\gamma',\gamma')\,dt\,.
\end{equation*}
Substituting now the Ricci soliton equation we get
\begin{equation*}
\int_0^L \sin^2\left({\frac{\pi t}{L}}\right) \frac{\mu}{n}\,\vert\gamma'\vert^2\,dt
-\int_0^L \sin^2\left({\frac{\pi t}{L}}\right)
\nabla^2_{\gamma',\gamma'}\widetilde{f}\, dt
\leq\frac{\pi^2(n-1)}{L^2}\int_0^L \cos^2{\left(\frac{\pi t}{L}\right)}\,dt
\end{equation*}
that is,
\begin{align*}
 \frac{\mu}{n}\int_0^L&\,\sin^2\left({\frac{\pi t}{L}}\right)\,dt
-\int_0^L \sin^2\left({\frac{\pi t}{L}}\right) \frac{d^2\,}{dt^2}[\widetilde{f}(\gamma(t))]\,dt\\
=&\,\frac{\mu}{n}\int_0^L \sin^2\left({\frac{\pi t}{L}}\right)\,dt
+2\frac{\pi^2}{L^2}\int_0^L \left[\sin^2\left({\frac{\pi t}{L}}\right)
-\cos^2\left({\frac{\pi t}{L}}\right)\right]\widetilde{f}(\gamma(t))\,dt\\
\leq &\,\frac{\pi^2(n-1)}{L^2}\int_0^L \cos^2\left({\frac{\pi t}{L}}\right)\,dt\,.
\end{align*}
As $\vert \widetilde{f}\vert\leq \max_M |f|\leq C$ and 
setting $A=\int_0^L \sin^2\left({\frac{\pi
    t}{L}}\right)\,dt=\int_0^L \cos^2\left({\frac{\pi t}{L}}\right)\,dt$ we obtain
\begin{equation*}
 \frac{\mu}{n} A - 2CA\frac{\pi^2}{L^2}\leq  \frac{(n-1)\pi^2}{L^2}A
\end{equation*}
hence,
\begin{equation*}
 L^2\leq\frac{n\pi^2(n-1+2C)}{\mu}\,.
\end{equation*}
This estimate says that all the counterimages of a point $p\in M$
belong to a bounded, hence compact, subset of $\widetilde{M}$. Since
such a set is discrete, it must be finite. The thesis then
follows. Also, the universal covering is compact.

Notice that, as a byproduct of this argument, we have that if the
potential function $f$ of a complete, shrinking gradient Ricci soliton
is bounded then the soliton is compact. Moreover, by
 equation~\eqref{equ8} it also follows that if $\RRR$ is bounded and
 $\vert\nabla f\vert$ is bounded, again the soliton is compact.
\end{proof}

\subsection{Another Proof of Proposition~\ref{p2}}
We give now a direct proof of Propositions~\ref{p2} and~\ref{p5} 
only using the defining equation~\eqref{soliton}, without passing by
Theorem~\ref{p4}.

We start with the following computation,
\begin{align*}
\Delta\nabla_i\omega_j=&\,\nabla_k\nabla^k\nabla_i\omega_j
=\nabla_k(\nabla_i\nabla^k\omega_j+\RRR^k_{\phantom{k}ijs}\omega^s)\\
=&\,-\nabla_k\nabla_i\nabla_j\omega^k-2\nabla_k\nabla_i\RRR^k_j+
\nabla_k\RRR^k_{\phantom{k}ijs}\omega^s+\RRR^k_{\phantom{k}ijs}\nabla_k\omega^s\\
=&\,-\nabla_i\nabla_k\nabla_j\omega^k-\RRR_{kijs}\nabla^s\omega^k
-\RRR_{ki\phantom{k}s}^{\phantom{ki}k}\nabla_j\omega^s
-2\nabla_i\nabla_k\RRR^k_j-2\RRR_{ki\phantom{k}s}^{\phantom{ki}k}\RRR^s_j
-2\RRR_{kij}^{\phantom{kij}s}\RRR^k_s\\
&\,-\nabla_j\RRR_{sk\phantom{k}i}^{\phantom{sk}k}\omega^s
-\nabla_s\RRR_{kj\phantom{k}i}^{\phantom{sk}k}\omega^s
+\RRR^k_{\phantom{k}ijs}\nabla_k\omega^s\\
=&\,-\nabla_i\nabla_j\nabla_k\omega^k-\nabla_i\RRR_{js}\omega^s
-\RRR_{js}\nabla_i\omega^s+\RRR_{ikjs}\nabla^s\omega^k
-\RRR_{is}\nabla_j\omega^s
-\nabla_i\nabla_j\RRR\\
&\,-2\SSS_{ij}
+2\RRR_{ikj}^{\phantom{ikj}s}\RRR^k_s+\nabla_j\RRR_{is}\omega^s
-\nabla_s\RRR_{ij}\omega^s
-\RRR_{i\phantom{k}js}^{\phantom{i}k}\nabla_k\omega^s\\
=&\,\nabla_i\nabla_j\RRR
-\nabla_i\RRR_{js}\omega^s
-\RRR_{js}\nabla_i\omega^s
+\RRR_{ikjs}\nabla^s\omega^k
-\RRR_{is}\nabla_j\omega^s
-\nabla_i\nabla_j\RRR\\
&\,-2\SSS_{ij}
+2\RRR_{ikj}^{\phantom{ikj}s}\RRR^k_s+\nabla_j\RRR_{is}\omega^s
-\nabla_s\RRR_{ij}\omega^s
-\RRR_{i\phantom{k}js}^{\phantom{i}k}\nabla_k\omega^s\\
=&\,-\nabla_i\RRR_{js}\omega^s
-\RRR_{js}\nabla_i\omega^s
+\RRR_{ikjs}\nabla^s\omega^k
-\RRR_{is}\nabla_j\omega^s\\
&\,-2\SSS_{ij}
+2\RRR_{ikj}^{\phantom{ikj}s}\RRR^k_s+\nabla_j\RRR_{is}\omega^s
-\nabla_s\RRR_{ij}\omega^s
-\RRR_{i\phantom{k}js}^{\phantom{i}k}\nabla_k\omega^s\\
\end{align*}
Then, we are ready to write the Laplacian of the Ricci tensor
\begin{align*}
2\Delta\RRR_{ij}=&\,-\Delta\left({\nabla_i\omega_j+\nabla_j\omega_i}\right)\\
=&\,\nabla_i\RRR_{js}\omega^s
+\RRR_{js}\nabla_i\omega^s
-\RRR_{ikjs}\nabla^s\omega^k
+\RRR_{is}\nabla_j\omega^s\\
&\,+\nabla_j\RRR_{is}\omega^s
+\RRR_{is}\nabla_j\omega^s
-\RRR_{jkis}\nabla^s\omega^k
+\RRR_{js}\nabla_i\omega^s\\
&\,+2\SSS_{ij}
-2\RRR_{ikj}^{\phantom{ikj}s}\RRR^k_s
-\nabla_j\RRR_{is}\omega^s
+\nabla_s\RRR_{ij}\omega^s
+\RRR_{i\phantom{k}js}^{\phantom{i}k}\nabla_k\omega^s\\
&\,+2\SSS_{ji}
-2\RRR_{jki}^{\phantom{jki}s}\RRR^k_s
-\nabla_i\RRR_{js}\omega^s
+\nabla_s\RRR_{ji}\omega^s
+\RRR_{j\phantom{k}is}^{\phantom{j}k}\nabla_k\omega^s\\
=&\,2\nabla_s\RRR_{ij}\omega^s+2\RRR_{js}\nabla_i\omega^s+2\RRR_{is}\nabla_j\omega^s\\
&\,-\RRR_{ikjs}\nabla^s\omega^k -\RRR_{jkis}\nabla^s\omega^k
+\RRR_{ikjs}\nabla^k\omega^s
+\RRR_{jkis}\nabla^k\omega^s\\
&\,+4\SSS_{ij}-4\RRR_{ikjs}\RRR^{ks}\,.\\
\end{align*}
Now, noticing that all the second line cancels,
\begin{align*}
\Delta\RRR_{ij}
=&\,\nabla_s\RRR_{ij}\omega^s
+\RRR_{js}\nabla_i\omega^s
+\RRR_{is}\nabla_j\omega^s
+2\SSS_{ij}
-2\RRR_{ikjs}\RRR^{ks}\\
=&\,\left\langle\nabla\RRR_{ij}\,\vert\,\omega\right\rangle-2\RRR_{ikjs}\RRR^{ks}\\
&\,-\frac{\RRR_{is}}{2}\Bigl(\nabla_j\omega^s+\nabla^s\omega_j
-\frac{2\mu}{n}g_{sj}\Bigr)
-\frac{\RRR_{js}}{2}\Bigl(\nabla_i\omega^s+\nabla^s\omega_i
-\frac{2\mu}{n}g_{si}\Bigr)\\
&\,+\RRR_{js}\nabla_i\omega^s
+\RRR_{is}\nabla_j\omega^s\\
=&\,\left\langle\nabla\RRR_{ij}\,\vert\,\omega\right\rangle-2\RRR_{ikjs}\RRR^{ks}
+\frac{2\mu}{n}\RRR_{ij}+
\frac{\RRR_{is}}{2}\bigl(\nabla_j\omega^s-\nabla^s\omega_j\bigr)
+\frac{\RRR_{js}}{2}\bigl(\nabla_i\omega^s-\nabla^s\omega_i\bigr)\,.
\end{align*}
Finally, contracting this equation with $g^{ij}$ we get
\begin{align}\label{equ99}
\Delta\RRR
=&\,\left\langle\nabla\RRR\,\vert\,\omega\right\rangle-2\RRR_{ks}\RRR^{ks}
+\frac{2\mu}{n}\RRR
+\frac{\RRR_{is}}{2}\bigl(\nabla_j\omega^s-\nabla^s\omega_j\bigr)g^{ij}
+\frac{\RRR_{js}}{2}\bigl(\nabla_i\omega^s-\nabla^s\omega_i\bigr)g^{ij}\\
=&\,\left\langle\nabla\RRR\,\vert\,\omega\right\rangle+\frac{2\mu}{n}\RRR-2\SSS\nonumber
\end{align}
by the skew--symmetry of the sum of the last two terms.\\
When $n=3$ or in general if the Weyl tensor is zero, as before, setting $\lambda_{\min}:M\to\R$ to be 
the minimal eigenvalue of the Ricci tensor, if $p\in M$ is 
the point where $\lambda_{\min}/\RRR$ gets its minimum with
eigenvector $v_p$, we consider a
local unit smooth tangent vector field $w=w^i$ such that $w(p)=v_p$,
$\nabla w^i(p)=\Delta w^i(p)=0$. Then the smooth function
$\RRR_{ij}w^iw^j/\RRR$ has a local minimum at $p$, 
$(\RRR_{ij}w^iw^j/\RRR)(p)=\lambda_{\min}(p)/\RRR(p)$,
$\nabla(\RRR_{ij}w^iw^j/\RRR)(p)=0$ and 
$\Delta(\RRR_{ij}w^iw^j/\RRR)(p)\geq 0$.\\
By the assumptions on the derivatives of $w$ at $p$ we have$\nabla(\RRR_{ij}/\RRR)(p)v_p^iv_p^j=0$ and 
$\Delta(\RRR_{ij}/\RRR)(p)v_p^iv_p^j\geq 0$, hence, 
\begin{align*}
0\leq&\,\Delta(\RRR_{ij}/\RRR)(p)v_p^iv_p^j\\
=&\,\frac{2\left(\RRR^3-n\lambda_{\min}\RRR^2+2(n-1)\lambda_{\min}^2\RRR-(n-1)\SSS\RRR
+(n-1)(n-2)\lambda_{\min}\SSS\right)}{\RRR^2}\\
&\,+\frac{\RRR_{is}}{2}\bigl(\nabla_j\omega^s-\nabla^s\omega_j\bigr)v_p^iv_p^j
+\frac{\RRR_{js}}{2}\bigl(\nabla_i\omega^s-\nabla^s\omega_i\bigr)v_p^iv_p^j\\
=&\,\frac{2\left(\RRR^3-n\lambda_{\min}\RRR^2+2(n-1)\lambda_{\min}^2\RRR-(n-1)\SSS\RRR
+(n-1)(n-2)\lambda_{\min}\SSS\right)}{\RRR^2}\\
&\,+\frac{\lambda_{\min}}{2}\bigl(\nabla_j\omega_s-\nabla_s\omega_j\bigr)v^s_pv^j_p
+\frac{\lambda_{\min}}{2}\bigl(\nabla_i\omega_s-\nabla_s\omega_i\bigr)v^s_pv^i_p\,.
\end{align*}
Again, by skew--symmetry the last two terms cancel and we get
$$
0\leq 
\RRR^3-n\lambda_{\min}\RRR^2+2(n-1)\lambda_{\min}^2\RRR-(n-1)\SSS\RRR
+(n-1)(n-2)\lambda_{\min}\SSS\,.
$$
Since this inequality and equation~\eqref{equ99} are 
respectively analogous to~\eqref{equ9} and~\eqref{equ4} 
for {\em gradient} Ricci solitons, following the proof of
Proposition~\ref{p2} we get directly to the conclusion only supposing
that $(M,g)$ is a Ricci soliton, without using Theorem~\eqref{p4} to
know that we are actually dealing with a {\em gradient} Ricci soliton.

\appendix

\section{Minimizing Perelman's Functional}\label{appA}
\def\WWW{{\mathcal{W}}}

We suppose that $(M,g)$ is a connected manifold, otherwise we work on
every single connected component. Let $dV$ be the Riemannian measure on $M$
associated to $g$.

We show here the existence of a smooth function $f:M\to\R$ such that 
\begin{equation*}
\RRR+2\Delta f-\vert\nabla f\vert^2+2\mu f/n
\end{equation*}
is constant, for every value $\mu>0$.

This is related to minimizing the following Perelman's functional
(see~\cite{perel1}),
\begin{equation}\label{perew}
\WWW (g ,f;\tau ) = \int_M \left[\tau (|\nabla f|^2 + \RRR) + f -n
\right] (4\pi \tau )^{-n/2} e^{-f}\, dV\,,
\end{equation}
where $\tau >0$ is a scale parameter, under the constraint
$$
f\in\Bigl\{f\in C^{\infty} (M)\,\Bigr\vert\,\int_M (4\pi \tau
)^{-{n/2}} e^{-f}\,dV=1\Bigr\}\,.
$$
We consider then the functional (which differs only for a constant
term by Perelman's one after the change of variable $u=e^{-f/2}$,
setting $\tau = n/2\mu$, and multiplying by $1/\tau$)
$$
{\mathcal F}(u) = \int_M \Bigl(\RRR u^2+4|\nabla u|^2
  -\frac{2\mu}{n}u^2\log u^2\Bigr)\,dV\,,
$$
then we look for the following constrained infimum
$$
\sigma= \inf_{u\in{\mathcal U}}{\mathcal{F}}(u)
$$
where $u\in C^\infty(M)$ runs in 
$$
{\mathcal{U}}= \Bigl\{u\in C^{\infty} (M)\,\Bigr\vert\,\int_M
u^2\,dV=1\,\,\text{and $u>0$}\Bigr\}
$$
(notice that the function $x^2\log x^2$ belongs to $C^1(\R)$ so the
integrand is well defined for functions in ${\mathcal U}$).

\begin{prop}\label{firsteigen}
The infimum $\sigma$ is finite and there exists a smooth nonnegative 
function $u\in {\mathcal{U}}$ achieving it.
\end{prop}
\begin{proof}
We show that $\sigma>-\infty$ and that 
any minimizing sequence $\{ u_i\}$ must be uniformly bounded in the
$H^1(M)$--norm.\\
We observe that for any $u\in C^\infty(M)$, by applying
Jensen inequality with respect to the probability measure 
$u^2\,dV$, one has
\begin{equation*}
\begin{aligned}
\int_M u^2\log u^2\,dV &=\frac{n-2}{2}\int _M u^2\log u^{\frac{4}{n-2}}\,dV\\
&\leq \frac{n-2}{2}\log\left(\int_M u^2 u^{\frac{4}{n-2}}\,dV\right)\\
&= \frac{n-2}{2}\log\left(\int_M u^{\frac{2n}{n-2}}\,dV\right)\\
&=\frac{n-2}{2}\log \left(\int _M  u^{2^*}\,dV\right)\,.
\end{aligned}
\end{equation*}
On the other hand, one has
\begin{equation*}
\begin{aligned}
\log\left(\int _M  u^{2^*}\,dV\right)
&\leq\log\left[\left(C_M\int_M(|\nabla
    u|^2+u^2)\,dV\right)^{\frac{n}{n-2}}\right]\\
& =\frac{n}{n-2} \log \left(C_M \int_M ( |\nabla u|^2
  +u^2)\,dV\right)\,,
\end{aligned}
\end{equation*}
where $C_M$ is the Sobolev constant of $(M,g)$.\\
Putting together these two inequalities we get
\begin{equation*}
  -\frac{2\mu}{n}
\int_M u^2\log u^2 \,dV\geq -\frac{2\mu}{n}\log \left(C_M \int_M ( |\nabla u|^2
  +u^2)\,dV\right)\geq -\int_M ( |\nabla u|^2+u^2)\,dV -C\,,
\end{equation*}
for some positive constant $C$ (depending only on $(M,g)$ and $\mu$). Hence,
\begin{equation}\label{wbound}
\begin{aligned}
{\mathcal F}(u) =& \int_M \Bigl(\RRR u^2+4|\nabla u|^2
  -\frac{2\mu}{n}u^2\log u^2\Bigr)\,dV\\
\geq&\,4\int _M (|\nabla u|^2+u^2)\,dV
+\int_M(\RRR_{\mathrm{\min}}-4) u^2\,dV
-\int_M(|\nabla u|^2 +u^2)\,dV-C\\
\geq&\,3\int _M |\nabla u|^2\,dV+(\RRR_{{\mathrm{\min}}}-1-C)\int_M u^2\,dV\\
\geq&\,\RRR_{{\mathrm{\min}}}-1-C\,,
\end{aligned}
\end{equation}
where in the last passage we used that $\int _M u^2\,dV=1$. This shows
that $\sigma>-\infty$.\\
The same argument gives that if $u_i\in C^\infty(M)$ is a minimizing
sequence for ${\mathcal {F}}$ such that $\Vert u\Vert_{L^2}=1$, then
$u_i$ is bounded in $H^1(M)$. 
Hence, we can extract a subsequence (not relabeled) weakly converging
in $H^1(M)$ and strongly converging in $L^{2+\varepsilon}(M)$, for some
$\varepsilon>0$, to some function $u$. Clearly, by the 
the $L^2$ convergence, we have $\int_M u^2\,dV=1$ and we can also
assume $u\geq 0$, by the definition of ${\mathcal {F}}$.\\
It is easy to see that the functional ${\mathcal {F}}$ is lower
semicontinuous with respect to the weak convergence in $H^1(M)$, as
the term $u^2\log u^2$ is 
{\em subcritical} (and the function $x^2\log x^2$ is continuous) hence
its integral is continuous.\\
Then, a limit function $u:M\to\R$ is a nonnegative, 
constrained minimizer of ${\mathcal F}$ in $H^1(M)$.\\
The Euler--Lagrange equation for $u$ read
\begin{equation*}
-4\Delta u + \RRR u -\frac{2\mu}{n}(u\log u^2 + u) =  Cu\,,
\end{equation*}
for some constant $C$. It can be rewritten as
\begin{equation}\label{euler-lagrange}
\Delta u = \RRR u/4 + Cu - \frac{\mu}{n}u\log u\,,
\end{equation}
to be intended in $H^1(M)$.\\
As $u$ is in $H^1(M)$ and the term $u^2\log u$ is subcritical, a
bootstrap argument together with standard elliptic estimates gives
that $u\in C^{1,\alpha}$.

Rothaus proved in~\cite{rothaus} that a solution to 
equation~\eqref{euler-lagrange} is positive or identically zero 
(see Appendix~\ref{appB}), 
this second possibility is obviously excluded by the constraint
$\int_M u^2\,dV=1$.\\ 
Hence, as $x^2\log x^2$ is smooth in 
$\R\setminus\{0\}$ we can infer that the function $u$ is
actually smooth.
\end{proof}

We can consider then the smooth function $f=-2\log u$. A simple
computation, using equation~\eqref{euler-lagrange} gives the
following.
\begin{cor}
For every $\mu>0$, there exists a smooth function $f:M\to\R$ such that
\begin{equation*}
\RRR+2\Delta f-\vert\nabla f\vert^2+2\mu f/n
\end{equation*}
is constant.
\end{cor}
Clearly, this argument works with every positive $\tau$ in the
functional $\WWW$.
\begin{cor}
For every $(M,g)$ and $\tau>0$ 
there exists a smooth function $f:M\to\R$ minimizing Perelman's
functional~\eqref{perew}.
\end{cor}

\section{Strong Maximum Principle for Semilinear 
Equations on Manifolds}\label{appB}

We consider the following elliptic semilinear equation
on a Riemannian manifold $(M,g)$
\begin{equation}\label{probsmp}
\Delta u(q)=\varphi(u(q),q)
\end{equation}
where $\varphi:\R\times M\to\R$ is a continuous function such that
$\varphi(0,q)=0$ for every $q\in M$.

If $u:M\to\R$ is a nonnegative $C^{1,\alpha}$ 
solution and $\varphi$ is $C^1$
then, in every connected component of $M$, 
either $u>0$ or $u$ is identically zero, by the strong maximum
principle.\\
If the function $\varphi$ is only continuous this is not true, as one
can see considering the $C^3$ function $f:\R\to\R$
\begin{equation*}
\begin{cases}
\begin{array}{lll}
f(x)=0\qquad &&\text{if $x<0$}\\
f(x)=x^4\qquad &&\text{if $x\geq0$}
\end{array}
\end{cases}
\end{equation*}
satisfying the equation $f^{''}=12\sqrt{f}$. Here
$\varphi(t)=12\sqrt{\vert t\vert}$ which is H\"older continuous but not
$C^1$.

Inspired by an analogous condition for uniqueness in ODE's, we have the
following proposition.
\begin{prop}\label{SMXODE}
If $\varphi\in C^0(\R\times M)$ and for every $p\in M$ there exist 
some $\delta>0$ and a continuous, nonnegative, concave function 
$\widetilde{\varphi}:[0,\delta)\to\R$ such that 
$\widetilde{\varphi}(0)=0$, $\widetilde{\varphi}(t)\geq\varphi(t,q)$
for every $t\in[0,\delta)$ and $q$ in some neighborhood of $p$, and
$\int_0^\delta\frac{dt}{\widetilde{\varphi}(t)}=+\infty$, then a 
nonnegative $C^{1,\alpha}$ solution $u$ of equation~\eqref{probsmp} is
either positive or identically zero in every connected component of $M$.
\end{prop}

\begin{rem} Notice that this condition is not very restrictive, in
  particular it holds for all the equations $\Delta u=\varphi(u)$ with
  $\varphi:\R\to\R$ Lipschitz.
\end{rem}

\begin{proof}
We suppose that $u=0$ at some point $p\in M$. This clearly implies
that $p$ is a minimum point for $u$, hence $\nabla u(p)=0$.\\
Letting $B_r=B_r(p)$ be the
geodesic ball of radius $r>0$ around $p\in M$, we define, for $r>0$
small (less than the injectivity radius $R$ of $p$ and such that $B_r$
is inside the neighborhood of the hypothesis about the function
$\varphi$),
\begin{align*}
S(r)=&\,\int_{\partial B_r} d\omega=\int_{\SSS^{n-1}}J(\theta,r)\,d\theta\\
h(r)=&\,\frac{1}{S(r)}\int_{\partial B_r} u\,
d\omega=\frac{1}{S(r)}\int_{\SSS^{n-1}}u(\theta,r)J(\theta,r)\,d\theta
\end{align*}
where $d\omega$ is the induced measure on the geodesic sphere
$\partial B_r$, $\theta$ is the coordinate on $\SSS^{n-1}$ with
canonical measure $d\theta$ and $J(\theta,r)$ is the density of
$d\omega$ with respect to $d\theta$.\\
By the assumptions on $u$, the function $h:[0,R)\to\R$ is $C^1$ and
nonnegative, moreover, since $u(p)=0$ and $\nabla u(p)=0$ we have
$h(0)=h^\prime(0)=0$.\\
By a standard computation, using the divergence theorem,
\begin{align*}
{h'(r)}=&\,\frac{\int_{\partial B_r} \frac{\partial u}{\partial r}\,d\omega
+\int_{\SSS^{n-1}}u(\theta,r)\frac{\partial J(\theta,r)}{\partial
  r}\,d\theta-h(r){S'(r)}}{S(r)}\\
=&\,\frac{\int_{B_r} \Delta u\,dV
+\int_{\SSS^{n-1}}u \frac{\partial\log{J(\theta,r)}}{\partial
  r}\,d\omega}{S(r)}-h(r)\frac{\partial \log{S(r)}}{\partial
  r}\\
\leq &\,\frac{\int_{B_r} \varphi(u(q),q)\,dV(q)}{S(r)}+
h(r)\left(\max_{\theta\in S^{n-1}}\frac{\partial\log{J(\theta,r)}}{\partial
  r}-\frac{\partial \log{S(r)}}{\partial
  r}\right)\\
\leq &\,\frac{\int_{B_r} \widetilde{\varphi}(u)\,dV}{S(r)}+
h(r)\max_{\theta\in
    S^{n-1}}\left\vert\frac{\partial\log({J(\theta,r)}/S(r))}{\partial
    r}\right\vert\\
= &\,\frac{\int_0^r \int_{\partial B_t} \widetilde{\varphi}(u)\,d\omega\,dt}{S(r)}+
h(r)\max_{\theta\in
    S^{n-1}}\left\vert\frac{\partial\log({J(\theta,r)}/S(r))}{\partial
    r}\right\vert\,.
\end{align*}

\begin{lemma}\label{geodiff} 
For every point $p$ in a Riemannian manifold $(M,g)$ we have 
\begin{equation*}
\max_{\theta\in
    S^{n-1}}\left\vert\frac{\partial\log({J(\theta,r)}/S(r))}{\partial
    r}\right\vert= o(1)
\end{equation*}
as $r$ tends to zero.
\end{lemma}
By means of this lemma we get
\begin{equation*}
{h'(r)}
\leq \frac{\int_0^r \int_{\partial B_t} \widetilde{\varphi}(u)\,d\omega\,dt}{S(r)}+Ch(r)
\end{equation*}
Using now Jensen inequality, by the concavity of $\widetilde{\varphi}$ on
$[0,\delta)$, for every $r>0$ small enough ($u(q)$ tends to zero as
$q\to p$) we get the following differential
inequality
\begin{equation*}
{h'(r)}\leq\frac{\int_0^r
  S(t)\widetilde{\varphi}(h(t))\,dt}{S(r)}+Ch(r)\,.
\end{equation*}
We introduce now the function
$\overline{h}(r)=\sup_{t\in[0,r]}h(t)$. Notice that $\overline{h}$ is
a nondecreasing continuous function, actually Lipschitz as $h\in
C^1$, and $\overline{h}\geq h$, moreover, at every differentiability point
(almost all by the Lipschitz property) we have either
$\overline{h}^\prime=0$ or $\overline{h}^\prime=h^\prime$.\\
We have then (in distributional sense),
\begin{align*}
{\overline{h}'(r)}\leq&\,\frac{\int_0^r
  S(t)\widetilde{\varphi}(h(t))\,dt}{S(r)}+Ch(r)\\
\leq&\,\frac{\int_0^r
  S(t)\widetilde{\varphi}(\overline{h}(t))\,dt}{S(r)}+C\overline{h}(r)\\
\leq&\,\frac{\widetilde{\varphi}(\overline{h}(r))\int_0^r
  S(t)\,dt}{S(r)}+C\overline{h}(r)\\
=&\,\frac{{\mathrm{Vol}}(B_r)}{S(r)}\widetilde{\varphi}(\overline{h}(r))
+C\overline{h}(r)\,.
\end{align*}
It is a standard fact that for $r>0$ small enough
${\mathrm{Vol}}(B_r)/{S(r)}\leq Cr\leq C$, hence we conclude
\begin{equation*}
\overline{h}'(r)\leq C\widetilde{\varphi}(\overline{h}(r))+C\overline{h}(r)\,,
\end{equation*}
for $r>0$ small enough.\\
It is now well known (by ODE's theory) that, since the condition 
$\int_0^\delta \frac{dt}{\widetilde{\varphi}(t)}=+\infty$
implies $\int_0^\delta
\frac{dt}{\widetilde{\varphi}(t)+t}=+\infty$, the function
$\overline{h}$ is identically 
zero on some interval $[0,\varepsilon)$.\\
The argument is easy: 
if $\overline{h}$ is positive for some interval $(0,\varepsilon)$ we
have (distributionally)
\begin{equation*}
\frac{\overline{h}^\prime}{\widetilde{\varphi}(\overline{h})+\overline{h}}\leq C
\end{equation*}
and integrating both sides,
\begin{equation*}
C\varepsilon\geq
\int_0^\varepsilon\frac{\overline{h}^\prime(r)}{\widetilde{\varphi}(\overline{h}(r))
+\overline{h}(r)}\,dr=\int_0^{\overline{h}(\varepsilon)}\frac{dt}{\widetilde{\varphi}(t)+t}
=+\infty
\end{equation*}
which is a contradiction.\\
This means that there exists a small $\varepsilon>0$ such that
$\overline{h}(r)=0$ on $[0,\varepsilon)$, hence, by construction, also
$h(r)=0$ on $[0,\varepsilon)$.

The same clearly holds for $u$ in a neighborhood of $p$, as $u\geq0$
and $h(r)$ is its spherical means on $\partial B_r(p)$. A standard
connectedness argument concludes the proof.
\end{proof}
\begin{proof}[Proof of Lemma~\ref{geodiff}]
We can compute the density of
the spherical measure $d\omega$ with respect to $d\theta$ as follows
(see~\cite[Sections~3.96,~3.98]{gahula}),
$$
J(\theta,r)=\sqrt{\det(g(Y_i(r),Y_j(r)))}
$$
where $Y_i(t)$ are Jacobi fields (that is,
$Y^{''}+\RRR(\gamma',Y)\gamma'=0$) along the geodesic
$\gamma(r)=\exp_p r\theta$ from $p\in M$, 
satisfying $Y(0)=0$ and $Y^\prime(0)=E_i$ where
$\{\theta,E_1,\dots,E_{n-1}\}$ is an orthonormal basis of $T_pM$.\\
By a standard computation, setting $G_{ij}(r)=g(Y_i(r),Y_j(r))$, we have
\begin{equation*}
\frac{\partial\log J(\theta,r)}{\partial r}
=\frac{\partial_rJ(\theta,r)}{J(\theta,r)}
=\frac{1}{2}\tr(\partial_rG(r)\circ G^{-1}(r))
=\tr(g(Y_i(r),Y_j^\prime(r))\circ G^{-1}(r))\,.
\end{equation*}
Now we expand the fields $Y_i$ and $Y^\prime_i$ by Taylor formula, around
$r=0$,
\begin{equation*}
\begin{cases}
Y_i(r)=rE_i(r) + o(r^2)\\
Y_i^\prime(r)= E_i(r)+o(r)
\end{cases}
\end{equation*}
where $E_i$ are the parallel fields along $\gamma(r)$ with
$E_i(0)=E_i$ and we used the formula
$Y^{''}+\RRR(\gamma',Y)\gamma'=0$. Moreover, the ``$o$--terms'' above
can be chosen independent of $\theta$.\\
Then, as the parallel transport is an isometry, we have
\begin{align*}
\frac{\partial\log J(\theta,r)}{\partial r}
=&\,\tr(g(rE_i(r) + o(r^2), E_j(r) + o(r))\circ G^{-1}(r))\\
=&\,\tr([r\delta_{ij}+o(r^2)]\circ [r^2\delta_{ij}+o(r^3)]^{-1})\\
=&\,r^{-1}\tr([\delta_{ij}+o(r)]\circ[\delta_{ij}+o(r)]^{-1})\,.
\end{align*}
Now, it is easy to see that
$[\delta_{ij}+o(r)]^{-1}=[\delta_{ij}+o(r)]$, hence
\begin{align*}
\frac{\partial\log J(\theta,r)}{\partial r}
=&\,r^{-1}\tr([\delta_{ij}+o(r)]\circ[\delta_{ij}+o(r)])\\
=&\,r^{-1}\tr[\delta_{ij}+o(r)]\\
=&\,\frac{n-1}{r}+o(1)\,.
\end{align*}
To conclude the proof of the lemma we notice that
\begin{align*}
\frac{\partial \log S(r)}{\partial r}
=\frac{\partial_r S(r)}{S(r)}
=&\,\frac{\partial_r\int_{{\mathbb
      S}^{n-1}}J(\theta,r)\,d\theta}{S(r)}\\
=&\,\frac{\int_{{\mathbb
      S}^{n-1}}\partial_rJ(\theta,r)\,d\theta}{S(r)}\\
=&\,\frac{\int_{{\mathbb
      S}^{n-1}}\frac{\partial_rJ(\theta,r)}{J(\theta,r)}J(\theta,r)
\,d\theta}{S(r)}\\
=&\,\frac{\int_{{\partial B_r}}
\frac{\partial_rJ(\theta,r)}{J(\theta,r)}\,d\omega}{S(r)}\,,
\end{align*}
hence, $\frac{\partial \log S(r)}{\partial r}$ is the spherical mean of
the function $\frac{\partial\log J(\theta,r)}{\partial r}=(n-1)/r +
o(1)$.\\
This finally implies that
\begin{align*}
\left\vert\frac{\partial\log({J(\theta,r)}/S(r))}{\partial
    r}\right\vert=&\,\left\vert\frac{\partial\log{J(\theta,r)}}{\partial
    r}-\frac{\partial\log{S(r)}}{\partial r}\right\vert\\
=&\,\vert(n-1)/r + o(1) -(n-1)/r + o(1)\vert\\
=&\, o(1)\,.
\end{align*}
\end{proof}

\begin{cor}\label{strper}
A nonnegative $H^1$ solution of the equation 
$$
\Delta u = \RRR u/4 + Cu - \frac{\mu}{n}u\log u\,,
$$
on a Riemannian manifold $(M,g)$ is either positive or identically
zero in every connected component of $M$.
\end{cor}
\begin{proof}
We already know that $u\in C^{1,\alpha}(M)$. In a neighborhood of any
$p\in M$ the scalar curvature $\RRR$
is bounded, then, $\Delta u(q)=\varphi(u(q),q)$ and locally 
$$
\varphi(t,q)=\RRR(q) t/4 + Ct - \frac{\mu}{n}t\log t\leq
\widetilde{C}t-\frac{\mu}{n}t\log t=\widetilde{\varphi}(t)\,,
$$
which is a continuous, concave, nonnegative function in some interval
$[0,\delta)$, satisfying the nonintegrability hypothesis in
Proposition~\ref{SMXODE}, by direct check. 
\end{proof}

\section{Some Open Problems}

\begin{prob} Is it possible to prove Theorem~\ref{p4} showing directly that the form
  $\omega$ is exact?
\end{prob}

\begin{prob}
Are there special conditions in dimension $n=4$ (on the Weyl tensor?)
assuring that a contracting, compact, Ricci soliton is trivial?
\end{prob}

\begin{prob}
Are there counterexamples in dimension $n=5$?
\end{prob}

\begin{prob}
The sectional curvatures of a compact, NONtrivial, contracting Ricci soliton 
can be all positive (nonnegative)?
\end{prob}

\begin{prob}
What are in general the properties of compact, NONtrivial, 
contracting Ricci solitons?
\end{prob}
See Derdzinski~\cite{derdz1}, Fern\'andez--L\'opez and
Garc\'{\i}a--R\'{\i}o~\cite{fergarciario}.

\bibliographystyle{amsplain}
\bibliography{Ricci}

\end{document}